%% file: courbe.tex
\documentclass[11pt]{article}
\input xy
\xyoption{all}
\usepackage{amsmath}
\usepackage{amsfonts}
\usepackage{amssymb}
\usepackage{stmaryrd}
\usepackage{graphicx,epsfig}
\usepackage{bbm,color}
\usepackage{amsthm}
\newtheorem{thm}{Theorem}[section]
\newtheorem{defn}[thm]{Definition}
\newtheorem{lem}[thm]{Lemma}

\newtheorem{prop}[thm]{Proposition}
\newtheorem*{cor*}{Corollary}
\newtheorem{rem}[thm]{Remark}
\newtheorem*{thm*}{Theorem}

\newcommand{\boS}{\mathcal{S}}
\newcommand{\boH}{\mathcal{H}}
\newcommand{\Z}{\mathbb{Z}}
\newcommand{\C}{\mathbb{C}}
\newcommand{\R}{\mathbb{R}}

\newcommand{\disj}{\text{\small{$\amalg$}}}

\renewcommand{\epsilon}{\varepsilon}
\renewcommand{\phi}{\varphi}

\DeclareMathOperator{\id}{Id}
\DeclareMathOperator{\tr}{Tr}

\newcommand{\ba}[1]{\overline{#1}}

\begin{document}
\title{Some asymptotics of TQFT via skein theory}

\author{Julien March\'e and Majid Narimannejad}

\date{}
\maketitle
\vspace{-20pt}

\begin{abstract}
For each oriented surface $\Sigma$ of genus $g$ we study a limit of quantum representations of the mapping class group arising in TQFT derived from the Kauffman bracket. We determine that these representations converge in the Fell topology to the representation of the mapping class group on $\boH(\Sigma)$, the space of regular functions on the $SL(2,\C)$ representation variety with its hermitian structure coming from the symplectic structure of the $SU(2)$-representation variety. As a corollary, we give a new proof of the  asymptotic faithfulness of quantum representations.
\end{abstract}

\section{Introduction}

A topological quantum field theory in dimension 2+1 is an algebraic structure very close to topology: roughly speaking, it associates to each surface a finite dimensional vector space and to each cobordism a linear map between the vector spaces associated to the boundaries. Such theories have physical origins: they were introduced by Witten \cite{witten} in the eighties from Chern-Simons actions and generated very rich mathematical developments. There are various rigorous constructions coming from geometric quantization, quantum groups and many other areas. Unfortunately, such constructions remain complicated and it is hard to make concrete computations. 

In this paper, we preferred the approach of \cite{BHMV} which defines TQFT in a purely combinatorial way: using skein theory and the Kauffman bracket, the authors defined a family of hermitian TQFT $(V_p,\langle\cdot,\cdot\rangle_p)$ corresponding for $p=2r$ to a $SU(2)$-theory with level $r-2$. Despite the simple and very beautiful structure of these combinatorial TQFT, the connection with geometry is less clear than from other approaches. In this article, we show that the same kind of connections can be found in a simple and direct way. From the axioms, a TQFT generates for any closed surface $\Sigma$ a family of representations of the extended mapping class group of $\Sigma$. In some sense, these representations carry the main topological meaning of TQFT, hence we would like to link them with some geometrical representation. 
The basic idea for this comes from a general belief that when $p$ goes to infinity, things become classical, by which we mean geometrical : such a belief is based on the so-called semi-classical approximation. Hence we propose to study the limit of $\rho_{2r}$, the quantum representations of $\Gamma_g$ on $V_{2r}(\Sigma)$ .

For this purpose, let us describe two classical spaces on which the mapping class group acts.

Fix a closed oriented surface $\Sigma$ of genus $g$. We call an isotopy class of 1-dimensional submanifold of $\Sigma$ a {\it multicurve}. The mapping class group of $\Gamma_g$ acts on the set of multicurves in a natural way. Call $C(\Sigma)$ the $\C$-vector space generated by multicurves: we obtain a representation of $\Gamma_g$ on $C(\Sigma)$. This fundamental representation carries almost all information about the structure of $\Gamma_g$. For instance, no non-trivial element of $\Gamma_g$ acts trivially on multicurves, except for the elliptic and hyperelliptic involutions in genus $1$ and $2$.

Another very natural spaces on which the mapping class group acts are the representation spaces of $\pi_1(\Sigma)$ in a fixed Lie group $G$. We note that $\boS(\Sigma,G)$ is isomorphic to $\hom(\pi_1(\Sigma),G)/G$. We are interested here in the cases $G=SU(2)$ and $G=SL(2,\C)$. These spaces have a rich structure: we will use the natural symplectic structure $\omega$ on the smooth part of $\boS(\Sigma,SU(2))$ and the  structure of an algebraic variety on $\boS(\Sigma,SL(2,\C))$. We define $\boH(\Sigma)$ as the ring of regular functions on $\boS(\Sigma,SL(2,\C))$. 

Using the natural inclusion of $\boS(\Sigma,SU(2))$ in $\boS(\Sigma,SL(2,\C))$, we can define a hermitian form on $\boH(\Sigma)$ by the formula

$$\langle f,g\rangle= \int_{\boS(\Sigma,SU(2))} f\ba{g} dV$$
Here, $dV$ is the volume form on $\boS(\Sigma,SU(2))$ induced by the symplectic form $\omega$.

We obtain the following result:
\begin{thm*}
Let $\Sigma$ be a closed oriented surface of genus $g$. For all even integers $p=2r$, there is a $\Gamma_g$-equivariant map $\psi_p:\boH(\Sigma)\to V_p(\Sigma)\otimes V_p(\Sigma)^*$ such that

One has $\langle v,w\rangle=\lim\limits_{p\to\infty}\frac{1}{r^{d(g)}}\langle\phi_p(v),\phi_p(w)\rangle_p$ for all $v,w\in\boH(\Sigma)$. Here, we have set $d(1)=1$ and $d(g)=3g-3$ for $g>1$.
This implies in particular that the quantum representations $\rho_{p}\otimes\ba{\rho_{p}}$ converge in the Fell topology to 
 $\rho:\Gamma_g\to U(\boH(\Sigma))$, the natural representation coming from the action of $\Gamma_g$ on $\boS(\Sigma,SL(2,\C))$.
\end{thm*}

As a corollary, we obtain a new proof of the result of \cite{freedman1} and \cite{andersen} about asymptotic faithfulness of quantum representations.

\begin{cor*}
Let $\Sigma$ be a closed oriented surface of genus $g$.
For any non-trivial $h$ in $\Gamma_g$ which is not the elliptic $(g=1)$ or hyperelliptic $(g=2)$ involution, there is some even $p_0$ such that $\rho_{p}(h)$ is not the identity for even $p>p_0$.
\end{cor*}

\begin{proof}
One can associate to any curve $\gamma$ on $\Sigma$ a regular function $f_{\gamma}$ on $\boS(\Sigma,SL(2,\C))$ by the formula $f_{\gamma}(\rho)=-\tr \rho(\gamma)$. For a disjoint union of curves, we associate the product of the functions associated to each component. In this way, we construct a map $f$ from $C(\Sigma)$ to $\boH(\Sigma)$.
By a result of \cite{bul} and \cite{prz}, the map $f$ is an isomorphism of vector spaces.  Therefore, we can think of a regular function on $\boS(\Sigma,SL(2,\C))$ as a linear combination of multicurves. 

Recall that no element of $\Gamma_g$ act trivially on $C(\Sigma)$ except the identity and the elliptic and hyperelliptic involutions in genus 1 and 2. Hence, we can suppose that there is some $v$ in $C(\Sigma)\simeq \boH(\Sigma)$ such that $w=hv-v$ is non-zero. This implies that $\langle w,w\rangle$ is non-zero, that is, the form $\langle\cdot,\cdot\rangle$ is non-degenerate.

In fact, if $\langle w,w\rangle=0$, the regular function on $\boS(\Sigma,SL(2,\C))$ associated to $w$ verifies $$\int_{\boS(\Sigma,SU(2))} |w|^2=0.$$ As $w$ is continuous, it must vanish on $\boS(\Sigma,SU(2))$. Moreover, as it is holomorphic on the space $\boS(\Sigma,SL(2,\C))$ and 0 on $\boS(\Sigma,SU(2))$, it vanishes identically. (See proof of theorem 1.4.1 in \cite{goldman2}).

Due to the equality $\langle w,w\rangle=\lim\limits_{r\to \infty}\frac{1}{r^{d(g)}}
\langle \phi_{2r} w ,\phi_{2r} w\rangle$, we can find $r_0$ such that for all $r\ge r_0$, $\phi_{2r} w \ne 0$. Hence $\phi_{2r}(hv)\ne \phi_{2r}(v)$ and $\rho_{2r}(h)$ cannot be the identity.
\end{proof}

\subsection{Plan of the proof of the theorem}

The heart of the proof is the construction of the map $\phi_p$, which is almost obvious, but is fundamental. As the space $\boH(\Sigma)$ is isomorphic to $C(\Sigma)$, to define a map $\phi_p$, it is sufficient to construct $\phi_p(\gamma)\in V_p(\Sigma)\otimes V_p(\Sigma)^*$ for any multicurve $\gamma$.

For such a multicurve, we consider the cobordism $\Sigma\times [0,1]$ with the multicurve embedded as $\gamma\times\{\frac{1}{2}\}$. The TQFT naturally induces an element $Z_p(\Sigma\times [0,1],\gamma)$ in 
$V_p(\Sigma \disj - \Sigma)=V_p(\Sigma)\otimes V_p(\Sigma)^*$. We call this element $\phi_p(\gamma)$.  This gives our fundamental map $\phi_p$, which is clearly equivariant because of the naturality of the construction.

To prove the theorem, one has to compute the limit of the expression $\frac{1}{r^{d(g)}}\langle \phi_p(\gamma),\phi_p(\delta) \rangle_p$ for two multicurves $\gamma$ and $\delta$.

We do this in two steps. In the first step, we assume that $\delta$ is empty. Using combinatorial techniques from \cite{BHMV}, we obtain for $\langle \phi_p(\gamma),1 \rangle_p$ an explicit fomula ressembling a Riemann sum. When we normalize it, it converges to an integral over a subspace of $\R^{d(g)}$, which we note $\langle \gamma \rangle$. By linearity, we extend $\langle \cdot\rangle$ to a map from $C(\Sigma)$ to $\C$.

In the second step, we use the connection between the TQFT $V_p$ and the Kauffman skein module at $A=-e^{\frac{i\pi}{p}}$. We find easily that $\frac{1}{r^{d(g)}}\langle \phi_p(\gamma),\phi_p(\delta) \rangle_p$ converges to $\langle \gamma\cdot \delta\rangle$, where $\cdot$ is the multiplication induced on $C(\Sigma)$ by its identification with the Kauffman skein algebra of $\Sigma\times [0,1]$ at $A=-1$ (see \cite{prz}).

On the other hand, it is well-known that this multiplication on $C(\Sigma)$ is isomorphic to the natural multiplication on $\boH(\Sigma)$, the space of regular functions on $\boS(\Sigma,SL(2,\C))$ (see \cite{bul,prz}).

It remains to identify the linear form on $\boH(\Sigma)$ defined by $f_{\gamma}\mapsto \langle \gamma\rangle$. Suppose that $\gamma$ is a multicurve. We choose curves $C_i$ on $\Sigma$ which decompose the surface into pants such that all components of $\gamma$ are parallel to some $C_i$. It is well known that the maps $f_i=f_{C_i}$ form a system of Poisson commuting functions on $\boS(\Sigma,SU(2))$.

As shown in  \cite{jeffrey}, the maps $(f_i):\boS(\Sigma,SU(2))\to \R^{d(g)}$ are the moment map for an action of a torus of dimension $d(g)$. By the Duistermaat-Heckman theorem, we obtain an explicit formula for the volume form $dV$ on $\boS(\Sigma,SU(2))$ and obtain finally the following striking formula:

$$\langle \gamma\rangle = \int_{\boS(\Sigma,SU(2))} f_{\gamma} dV.$$

From this formula, we deduce the theorem.

\subsection{Remarks and perspectives}

The main motivation for this work came from the article \cite{freedman2} about the asymptotics of quantum representations of the mapping class group of the torus. Our approach is different in the sense that we study the limit of $V_p\otimes V_p^*$ instead of simply $V_p$. We were also inspired by the ideas contained in the paper \cite{magnet}. Our work is of course related to the article \cite{andersen} where similar ideas appear, and has also some intersection with \cite{yangmills}.

There are many questions naturally linked to our results:

\begin{itemize}
\item[-]
How can we link our asymptotic result to the asymptotics considered in \cite{freedman2}?
\item[-]
Can we apply our result or some refinements to the problem of \cite{masbaum}? The Nielsen-Thurston classification of the elements of the mapping class group is directly related to their action on multicurves. As the quantum representations converge to this action, can we find some trace of this classification in quantum representations?
\item[-]
In \cite{BHMV}, one can choose any primitive $4r$-root of unity to construct a TQFT. We have chosen roots converging to $-1$. Is it possible to develop the same asymptotics for roots of unity converging to different complex numbers?
\item[-]
Can we obtain a stronger convergence for the sequence involved in the theorem?
\end{itemize}

\section{Review of TQFT}

This part is a quick and formal review of TQFT constructed in \cite{BHMV} which we give to fix notations and settings,  and to recall results that will be used in this paper. We refer the interested reader to the beautiful original paper.

Fix an even integer $p=2r$. The complex number $A=-e^{i\pi/2r}$ is a primitive $4r$-th root of unity. One can construct from it a 2+1 topological quantum field theory.

In the notations of \cite{BHMV}, we set $\kappa=e^{-\frac{i\pi}{2r}-\frac{i\pi(2r+1)}{12}}$ and 
$\eta=\sqrt{\frac{2}{r}}\sin(\frac{\pi}{r})$. We define $C_r=\{0,1,\cdots,r-2\}$ which will be called the set of {\it colors}.

A triple $(a,b,c)$ of elements of $C_r$ is called {\it $r$-admissible} if $a+b+c$ is even, the triangle inequality $|a-b|\le c\le a+b$ is satisfied and moreover we have $a+b+c<2r-2$. 

For any integer $j$, we denote the quantum integer $[j]$ given by the formula $[j]=\frac{\sin(\frac{\pi j}{r})}{\sin(\frac{\pi}{r})}$, and define the quantum factorial by the formula $[j]!=\prod_{k=1}^j[k]$.

\subsection{The cobordism category}

A TQFT is a linear representation of a cobordism category. In our settings, the objects of our category are oriented surfaces with marked points and $p_1$-structures. 
\begin{itemize}
\item[-] A marking of a surface $\Sigma$ is a family $(z_j,c_j)_{j\in J}$ where $(z_j)$ is a family of distinct points in $\Sigma$ with for all $j\in J$ a non zero tangent direction $v_j$ at $z_j$ on $\Sigma$. For all $j\in J$, $c_j$ is a color in $C_p$.
\item[-] A $p_1$-structure is a somewhat complicated object, used to solve the so-called ``framing anomaly''. Consider the map $p_1:BO\to K(\Z,4)$ corresponding to the first Pontryagin class. Let $X$ be its homotopy fiber, i.e. the set of couples $(x,\gamma)\in BO\times C([0,1],K(\Z,4))$ satisfying $\gamma(0)=*$, and $\gamma (1)=p_1(x)$. Let $E$ be the universal stable bundle over $BO$, and $E_X$ its pull-back over $X$. A $p_1$-structure on a manifold $M$ is a fiber map from the stable tangent bundle of $M$ to $E_X$.
\end{itemize}
In the notation of an object $(\Sigma,z,c)$, we do not mention the directions $v_j$ and the $p_1$-structure, although they are present.

We define now morphisms : let $(\Sigma_1,z_1,c_1)$ and $(\Sigma_2,z_2,c_2)$ be two objects as defined above. A morphism is 
\begin{itemize}
\item[-]
An oriented 3-manifold $M$ whose boundary is decomposed as $\partial M = -\Sigma_1\disj \Sigma_2$, where $-\Sigma$ means $\Sigma$ with opposite orientation.
\item[-]
A colored banded trivalent graph $G$ embedded in $M$ whose restriction to the boundary is compatible with the marked points.
\item[-]
A $p_1$-structure on $M$ extending the $p_1$-structure given on the boundary.
\end{itemize}

A banded trivalent graph $G$ in $M$ is a 1-3-valent graph contained in an oriented surface $SG\subset M$ such that 
\begin{enumerate}
\item[(i)]
$G$ meets $\partial M$ transversally on the set of 1-valent vertices of $G$ noted $\partial G$.
\item[(ii)]
The surface $SG$ is a regular neighborhood of $G$ in $SG$, and $SG\cap\partial M$ is a regular neighborhood of $G\cap\partial M$ in $SG\cap\partial M$.
\end{enumerate}
A coloring of $G$ is a map $\sigma$ from the set of edges of $G$ to $C_r$ such that the colors of the edges meeting at each vertex are $r$-admissible. 
The restriction of a banded graph $G\subset M$ on some connected component of $\partial M$ gives marked points $(z_j)_{j\in J}$ with tangent directions $(v_j)_{j\in J}$, whereas the restriction of a coloring gives colors $(c_j)_{j\in J}$.

Two morphisms are called equivalent if the corresponding manifolds are isomorphic, the banded graphs are isotopic and the $p_1$-structures are homotopic relative to the boundary.

\subsection{Main properties of TQFT}


The theorem proved in \cite{BHMV} states that for each integer $p$, there is a functor $(V_p,Z_p)$ from the precedent cobordism category to the category of finite $\C$-vector spaces.

This means that to every object $(\Sigma,z,c)$ we can associate a vector space $V_p(\Sigma,z,c)$ and to any morphism $(M,G)$ between two objects a linear map $Z_p(M,G)$ between the vector spaces corresponding to the objects. By convention, $V_p(\emptyset)=\C$, hence any closed manifold $(M,G)$ acts as a scalar $\langle M,G\rangle_p$ which is a 3-manifold invariant. 
Moreover, there is natural hermitian form $\langle\cdot,\cdot\rangle_p$ on $V_p(\Sigma,z,c)$ such that for any two morphisms $(M_1,G_1)$ and $(M_2,G_2)$ from $\emptyset$ to $(\Sigma,z,c)$, we have 
$\langle Z_p(M_1,G_1),Z_p(M_2,G_2)\rangle_p=\langle M_1\cup(-M_2),G_1\cup G_2\rangle_p$.

We give here some important results related to this construction:

\begin{thm}[BHMV]\label{thmbhmv}
Let $(\Sigma,z,c)$ be a surface with marked points and $p_1$-structure. Let $H$ be a handlebody whose boundary is $\Sigma$ and with a $p_1$-structure extending that of $\Sigma$. Let $G$ be a 1-3-valent banded graph in $H$ such that $\partial G=z$ and such that $H$ is a tubular neighborhood of $G$. 
For each coloring $\sigma$ of $G$ compatible with the coloring of the boundary, we note $u_{\sigma}$ the element induced by $Z_p$ in $V_p(\Sigma,z,c)$. 

Then the elements $u_{\sigma}$ form an orthogonal basis of $V_p(\Sigma,z,c)$, and if $G$ does not contain any closed loop, we have 
$$
\langle u_{\sigma},u_{\sigma}\rangle_p =\eta^{\# v -\# e}
\frac{\prod_v \langle \sigma_v\rangle}{\prod_e \langle \sigma_e\rangle}$$
In this formula, $v$ ranges over the set of vertices of $G$ and $e$ over the set of $edges$. Moreover, for any trivalent vertex $v$, $\sigma_v$ is the triple of colors of the edges adjacent to this vertex and for any monovalent vertex $v$, $\sigma_v$ is the color of the edge incoming to it.

We set $\langle j\rangle=(-1)^j[j+1]$ and $\langle a,b,c\rangle=(-1)^{\alpha+\beta+\gamma}\frac{[\alpha+\beta+\gamma+1]![\alpha]![\beta]![\gamma]!}{[a]![b]![c]!}$ where $\alpha,\beta$ and $\gamma$ are defined by the equations $a=\beta+\gamma,b=\alpha+\gamma,c=\alpha+\beta$.

If $G$ is reduced to a closed loop, then the formula is simply $\langle u_\sigma,u_\sigma\rangle_p=1$.
\end{thm}
\begin{rem}
We check that for our choice of root of unity $A$ and for a surface $\Sigma$ without marked points, the hermitian pairing on $V_p(\Sigma)$ is definite positive.
\end{rem}

\subsection{Kauffman Bracket and TQFT}

We define $K(M)$ as the usual skein module of the manifold $M$. We refer to \cite{prz} for a complete account, but we will recall here what we need.
Let $A$ be some indeterminate. The $\Z[A,A^{-1}]$-module $K(M)$ is the free module generated by isotopy class of banded links in $M$ quotiented by the submodule generated by the local relations of figure 1.

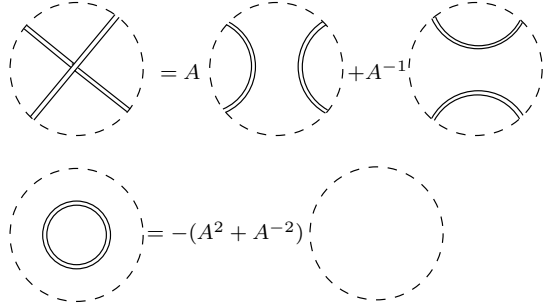
\begin{figure}[htbp]
\begin{center}
\input{kauf.pstex_t}
\caption{Kauffman relations}
\end{center}
\end{figure}

For any $u\in\C\setminus\{0\}$, we set $K(M,u)=K(M)\otimes_{\Z[A,A^{-1}]} \C$ where $A$ acts on $\C$ by multiplication by $u$.

The following proposition is a consequence of the construction of the TQFT.

\begin{prop}[Proposition 1.9 in \cite{BHMV}]\label{engendre}
Let $M$ be a connected $3$-manifold with $p_1$-structure and boundary $\Sigma$ (without boundary or marked points). Then there is a surjective map from $K(M,-e^{i\pi/p})$ to $V_p(\Sigma)$.

This map is defined by sending the element $L\otimes 1$ to $Z_p(M,L)$.
\end{prop}

\section{Convergence of TQFT}

\subsection{Settings}
Let $\Sigma$ be a closed oriented surface of genus $g$ with $p_1$-structure. We denote by $\Gamma_g$ the mapping class group of $\Sigma$. Fix $p=2r$.

If $h$ is an element of $\Gamma_g$, we can construct a cobordism $C_h$ from $\Sigma$ to itself as $\Sigma\times [0,1]$ where we identify the first boundary component with $\Sigma$ using the identity and the second one using $h$. 
If $h'$ is another element  of $\Gamma_g$, the cobordisms $C_h\circ C_{h'}$ and $C_{hh'}$ are diffeomorphic. We should obtain a  representation of $\Gamma_g$ on $V_r(\Sigma)$ by considering the linear map $Z_p(C_h)$. 
The problem is that we have not chosen any $p_1$-structure on $C_h$, and there is no canonical choice to make.

One way to get rid of this annoying fact is to consider the action of $\Gamma_g$ on $V_p(\Sigma)\otimes V_p(\Sigma)^*=V_p(\Sigma\disj-\Sigma)$. 
This action of $h$ on this space is given by $Z_p(C_h\disj -C_h)$ where we choose any $p_1$-structure on $C_h$ and put the same one on $-C_h$. The action does not depend any more on the $p_1$-structure: in fact, in a cobordism $M$, if we change the $p_1$-structure, the linear map $Z_p(M)$ is changed by a multiple of $\kappa$, a root of unity.  When we take the dual, the root becomes its conjugate. Hence, the two ``anomalies'' cancel and we get a true representation of $\Gamma_g$.

We thus obtain a sequence of representations $(V_p(\Sigma),Z_p)$ of $\Gamma_g$, and want to find their limit in some sense. The problem is that the spaces on which the mapping class group acts are a priori completely different. We need a way to compare them which is suggested by Proposition \ref{engendre}.

\begin{defn}
Let $\Sigma$ be a closed oriented surface. We call a 1-submanifold of $\Sigma$ without component bounding a disc  in $\Sigma$ a multicurve. 
We define $C(\Sigma)$ as the free $\C$-vector space generated by isotopy classes of multicurves on $\Sigma$. 
\end{defn}

Given a multicurve $\gamma$ in $\Sigma$, one can give it a banded structure by taking a neighborhood of it in $\Sigma$. We can consider the curve $\gamma$ as a banded link in $\Sigma\times [0,1]$ by sending it to $\gamma\times\{1/2\}$. We use the same notation for the multicurve on $\Sigma$ and its associated banded link in $\Sigma\times[0,1]$.

In \cite{prz}, it is shown that the Kauffman skein module $K(\Sigma\times [0,1])$ is a free $\Z[A,A^{-1}]$-module with basis the isotopy classes of multicurves. It provides an isomorphism of vector spaces between $C(\Sigma)$ and $K(\Sigma\times[0,1],u)$ for any $u$ in $\C\setminus\{0\}$. In particular, using Proposition \ref{engendre}, we get a surjective map 
$$\phi_p: C(\Sigma)\to K(\Sigma\times [0,1],-e^{i\pi/p})\to V_p(\Sigma\disj-\Sigma).$$

\begin{thm}\label{principal}
Let $\Sigma$ be a closed oriented surface of genus $g$. There is an hermitian pairing $\langle \cdot,\cdot\rangle$ on $C(\Sigma)$ such that for all $x$ and $y$ in $C(\Sigma)$, the following holds, where $d(1)=1$ and $d(g)=3g-3$ for $g>1$.

$$\langle x,y\rangle =\lim_{r\to \infty}\frac{1}{r^{d(g)}}\langle \phi_{p}(x),\phi_{p}(y)\rangle_{p}.$$
\end{thm}

\subsection{The trace function}

\begin{defn}
Let $\Sigma$ be a closed oriented surface of genus $g$ and $\gamma$ be a multicurve on $\Sigma$. 
We set $\tr_p(\gamma)=\langle \Sigma\times S^1,\gamma\rangle_p$. Here, $\gamma$ is seen as a banded link lying in the slice $\Sigma\times\{1\}$ of $\Sigma\times S^1$.
\end{defn}

\begin{lem}\label{trace}
Suppose that a surface $\Sigma$ is presented as the boundary of a handlebody $H$ which retracts on a trivalent banded graph $G$ as in Theorem \ref{thmbhmv}. We choose meridian disks $D_e$ transverse to each edge of $G$ and define $C_e=\partial D_e$ : the curves $C_e$ are disjoint on $\Sigma$. We choose a non-negative integer $m_e$ for each edge of $G$.

Then we define $\gamma$ as the multicurve on $\Sigma$ obtained by taking $m_e$ parallel copies of $C_e$ for each edge of $G$. We have 
$$\tr_p(\gamma)=\sum_{\sigma}\prod_e \left[-2\cos\left(\frac{(\sigma_e+1)\pi}{r}\right)\right]^{m_e}.$$
Here $\sigma$ ranges over $r$-admissible colorings of $G$ and $e$ ranges over edges of $G$.
\end{lem}

\begin{proof}
The proof is an easy consequence of the gluing axioms and the following fact from skein theory:  a trivial curve colored with 1 and making a Hopf link with a curve colored with $j$ may be removed and replaced by a factor $-A^{2j+2}-A^{-2j-2}=-2\cos(\frac{(j+1)\pi}{r})$. We refer for instance to Lemma 3.2 of \cite{BHMV1}.

For each edge $e$ of $G$, let us add a special component $S_e$ to $\gamma$ and cut $\Sigma$ along these special curves. The manifold $\Sigma\times S^1$ appears as a gluing of submanifolds isomorphic to a product of pants with $S^1$. Each gluing is realized by the trivial cobordism $S_e\times S^1\times [0,1]$ between two torus $S_e\times S^1$. 

For each boundary circle $S_e$, there is a preferred basis for $V_p(S_e\times S^1)$ given by $(D_e,i)\times S^1$ where $D_e$ is the disk bounding $S_e$ with one point marked with $i$. We call $e_i=Z_p((D_e,i)\times S^1)$ the corresponding basis element.
Thanks to the formulas of theorem \ref{thmbhmv}, this basis is orthonormal. Hence, the trivial cobordism satisfies $Z_p(S_e\times S^1\times [0,1])=\sum_i e_i\otimes e_i^*=\id$. This means that we can replace for each edge $e$  the trivial cobordism $S_e\times S^1\times [0,1]$ by a sum over $i$ of two solid tori $(D_e,i)\times S^1$ glued along the boundary components of the pants (times $S^1$) which have been cut. This is suggested in  the figure 2.

\begin{figure}[htbp]
   \centering
   \input{bretzel.pstex_t} 
   \caption{Contracting tensors in TQFT}
\end{figure}
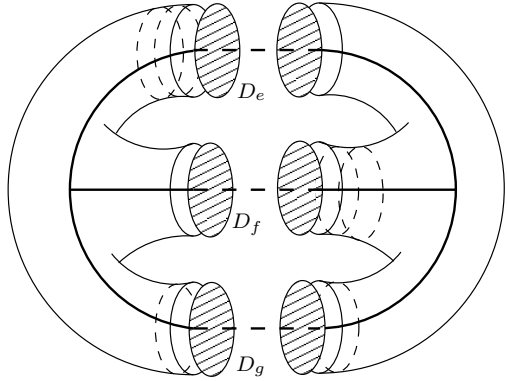

After having performed this decomposition along all the edges of $G$, we obtain a sum over all colorings $\sigma$ of the edges by elements of $C_r$. To each coloring is associated a disjoint union of manifolds isomorphic to $S^2\times S^1$ with colored links lying inside.

More precisely, the contribution of each pant with boundary 
$S_e\disj S_f\disj S_g$, $m_e$, $m_f$, $m_g$ parallel copies of the boundary components and disks colored respectively with $\sigma_e,\sigma_f$ and $\sigma_g$ may be explicitely computed.
The contribution is given by the $Z_p$ invariant of a sphere with 3 marked points with colors $\sigma(e),\sigma(f)$ and $\sigma(g)$ times a circle. The curves occuring as a product of the marked points with a circle are respectively linked to $m_e$, $m_f$ and $m_g$ trivial circles. As we have said in the beginning of the proof, these circles may be removed and replaced by the factor 
$$\left(-2\cos(\frac{(\sigma_e+1)\pi}{r})\right)^{m_e}\left(-2\cos(\frac{(\sigma_f+1)\pi}{r})\right)^{m_f}\left(-2\cos(\frac{(\sigma_g+1)\pi}{r})\right)^{m_g}.$$

What remains is $Z_p((S^2,\sigma_e,\sigma_f,\sigma_g)\times S^1)=\dim V_p(S^2,\sigma_e,\sigma_f,\sigma_g)$. But the latter space is either 1 dimensional if the triple $(\sigma_e,\sigma_f,\sigma_g)$ is $r$-admissible or is $\{0\}$. By taking the sum over all maps $\sigma$ from the set of vertices of $G$ to $C_r$, we finally obtain the formula of the lemma.
\end{proof}

\subsection{Limit of the trace function}

As before, fix a surface $\Sigma$ presented as in Theorem \ref{thmbhmv} as the boundary of a handlebody $H$ which retracts on a trivalent banded graph $G$.

The number of edges of $G$ is $3g-3$ if $g>1$ or $1$ if $g=1$. 
We denote this number by $d(g)$ and consider the subset $U_g$ of $\R^{d(g)}$ consisting of all maps $\tau$ from the set of edges of $G$ to $[0,1]$ such that for all triples of incoming edges $(e,f,g)$ of some vertex, the following relations are satisfied: 
\begin{itemize}
\item
$| \tau_f-\tau_g|\le \tau_e\le \tau_f+\tau_g$
\item
$\tau_e+\tau_f+\tau_g\le 2$
\end{itemize}

We use the formula of Lemma \ref{trace} to deduce the asymptotics of the trace function.
\begin{lem}\label{calcul}
With the  same hypothesis as Lemma \ref{trace}, let $F:U_g\to \R$ be the map defined by $F(\tau)=\prod_e(-2\cos(\tau_e\pi))^{m_e}$. Then

$$\lim_{r\to \infty} \frac{1}{r^{d(g)}}\tr_p(\gamma)=2^{g-d(g)}\int_{U_g}F(\tau)d\tau.$$

\end{lem}
\begin{proof}
The formula for $\tr_p(\gamma)$ looks like a Riemann sum, hence the result should not be a surprise. To obtain the precise result, we have to decompose $U_g$ into small pieces parametrized by $r$-admissible colorings $\sigma$. 

Given a positive integer $r$ and any coloring $\sigma$ from the set of edges of $G$ to $C_r$, we define the following set $A^r_{\sigma}=\prod_e [\frac{\sigma_e}{r},\frac{\sigma_e+1}{r})\subset \R^{d(g)}$. As $\sigma$ runs over $r$-admissible colorings of $G$, these sets do not cover $U_g$ because of the parity condition. We have to pack some sets $A^r_{\sigma}$ together, which we do in the following way.

Choose a subspace $S$ of $C_1(G,\Z_2)$ such that $C_1(G,\Z_2)=S\oplus Z_1(G,\Z_2)$. The subspace $S$ has dimension $d(g)-g$. For an admissible coloring $\sigma$ of $G$, we define $B^r_{\sigma}=\bigcup_{\rho \in S} A^r_{\sigma+\rho}$. Here we have identified $\Z/2\Z$ with the set $\{0,1\}$. It happens that the sets $B_{\sigma}^r$ are disjoint and cover $U_g$. We prove that they are disjoint: if we have $\sigma+\rho=\sigma'+\rho'$ with $\sigma$ and $\sigma'$ admissible and $\rho,\rho'$ in $S$, then consider these maps modulo 2. If we apply the boundary map, the admissible colorings vanish by definition, and we have $\partial \rho =\partial \rho'$. But $\partial$ induces a bijection from $S$ onto its image, hence we have $\rho=\rho'$, and it follows that $\sigma=\sigma'$. Hence the sets $B^r_{\sigma}$ are actually disjoint. Moreover the measure of $B^r_{\sigma}$ is $\frac{2^{d(g)-g}}{r^{d(g)}}$. It follows that $\sum\limits_{\sigma, r-\text{admissible}}F(\frac{\sigma_e+1}{r})\frac{2^{g-d(g)}}{r^{d(g)}}$ converges to $\int_{U_g}F(\tau)d\tau$ and the result is proved.
\end{proof}

\subsection{Proof of the theorem \ref{principal}}

Let $\Sigma$ be a closed oriented surface of genus $g$. 
We recall that $C(\Sigma)$ and $K(\Sigma\times[0,1],u)$ are isomorphic for any $u$ in $\C\setminus\{0\}$. The stacking product gives to $K(\Sigma\times[0,1])$ a natural algebra structure which induces an algebra structure on $C(\Sigma)$ for each $u\in\C\setminus\{0\}$. We consider the algebra structure obtained for $u=-1$. 

Fix $\gamma$ and $\delta$, two multicurves on $\Sigma$. We aim to compute the limit of $\frac{1}{r^{d(g)}}\langle \phi_r(\gamma),\phi_r(\delta)\rangle_r$ as $r$ goes to infinity. The right hand side is the quantum invariant of two thickened surfaces $\Sigma$ with a multicurve inside, glued along their boundary. Instead of gluing the two boundaries simultaneously, we glue one and then the other. If we glue one boundary component, we obtain the stacking product of $\gamma$ and $\delta$. 
In the skein module for generic $A$, we have a decomposition $\gamma\cdot \delta=\sum_i c_i \zeta_i$ for some multicurves $\zeta_i$ and some Laurent polynomials $c_i$ in  $\Z[A,A^{-1}]$.  When evaluating this combination in $V_p(\Sigma\disj -\Sigma)$, we have to specialize $A$ to $-e^{\frac{i\pi}{p}}$. 
In formulas, we have $\phi_p(\gamma\cdot \delta)=\sum_i c_i(-e^{\frac{i\pi}{p}}) \phi_p(\zeta_i)$. Then, we glue together the remaining boundary components and obtain $\langle \phi_p(\gamma),\phi_p(\delta)\rangle_p=\sum_i c_i(-e^{\frac{i\pi}{p}}) \tr_p(\zeta_i)$.

The asymptotic formula becomes clear if we define the following linear form on $C(\Sigma)$:

\begin{defn}\label{form}
Let $\gamma$ be a multicurve on $\Sigma$. Then there is a pants decomposition associated to $\gamma$ such that all components of $\gamma$ are parallel copies of the boundary circles.
As in lemma \ref{trace}, we define $\langle \gamma\rangle=2^{g-d(g)}\int_{U_g}F(\tau)d\tau$ where $F(\tau)=\prod_e(-2\cos(\tau_e \pi))^{m_e}$. The expression of $\langle\gamma\rangle$ as a limit shows that this definition does not depend on the pants decomposition. 
We extend $\langle\cdot\rangle$ to a linear form on $C(\Sigma)$.
\end{defn}

Coming back to our computation, we obtain:
$\lim\limits_{p\to\infty}\frac{1}{r^{d(g)}}\langle \phi_p(\gamma),\phi_p(\delta)\rangle_p= \sum_i c_i(-1)\langle \zeta_i\rangle=\langle \gamma\delta\rangle$. Finally, we define an hermitian form on $C(\Sigma)$ by the formula $\langle x,y\rangle= \langle x\ba{y}\rangle$ where the product corresponds to the skein module product for $A=-1$, and the conjugation corresponds to conjugation of coefficients in $C(\Sigma)$. We have proven the following result: 

For all $x,y\in C(\Sigma)$, we have $\langle x,y\rangle =\lim\limits_{p\to \infty}\frac{1}{r^{d(g)}}\langle \phi_p(x),\phi_p(y)\rangle_p.$

\section{Geometric interpretation}

The heart of the following geometric interpretation is the theorem of \cite{bul} and \cite{prz} stating that the algebra $K(\Sigma\times [0,1],-1)$ is isomorphic to $\boH(\Sigma)$, the ring of regular functions on the $SL(2,\C)$-representation variety of $\Sigma$. Recall that the isomorphism is given by $f_{\gamma}(\rho)=-\tr (\rho(\gamma))$ when $\gamma$ is a connected curve on $\Sigma$ and $\rho:\pi_1(\Sigma)\to SL(2,\C)$ is a representation of $\pi_1(\Sigma)$.

The space $C(\Sigma)$ is now identified, together with its algebra structure. It remains to identify the linear form $\langle\cdot\rangle$ of definition \ref{form}.

Recall that the $SL(2,\C)$-representation variety contains the $SU(2)$-representation variety, which carries a natural symplectic form $\omega$ defined in \cite{atiyah,goldman}. Following \cite{jeffrey}, we define 
$\boS_g$ to be the moduli space of irreducible representations of $\pi_1(\Sigma)$ on $SU(2)$ and $\ba{\boS}_g$ the moduli space of all representations. 
Then it is known that $\boS_g$ is a smooth $2d(g)$ manifold with symplectic form $\omega$ obtained by symplectic reduction from the form $\ba{\omega}(a,b)=\frac{1}{4\pi^2}\int_{\sigma}\tr a\wedge b$ for $a,b\in \Omega^1(\Sigma,su(2))$. 

We denote the volume form on $S_g$ by $dV=\frac{\omega^{d(g)}}{d(g)!}$.

\begin{prop}
For all multicurves $\gamma$ on $\Sigma$, we have
$$\langle \gamma \rangle = \int_{S_g} f_{\gamma} dV$$
\end{prop}
\begin{proof}
We give a proof of this proposition by adapting the results of \cite{jeffrey}. 

Fix a pants decomposition of $\Sigma$ associated to $\gamma$ and denote the set of curves bounding the pants by $C_e$. We define the function $h_e$ on $\ba{\boS}_g$ by the formula $\tr \rho(C_e)=2\cos(\pi h_e(\rho))$. 

Where the functions $h_e$ are not equal to 0 or 1, they Poisson commute and their Hamiltonian flows define a torus action on $\ba{\boS}_g$.

In fact, the function $h$ takes its values in $U_g$ and we have the following theorem:
\begin{thm}[\cite{jeffrey}]
Let  $U_g^{gen}$ be the interior of $U_g$ in $\R^{d(g)}$.

For $x$ in $\ba{\boS}_g$ such that $h(x)=y\in U_g^{gen}$, the torus action identifies $h^{-1}(y)$ with $U(1)^{d(g)}/\Z_2^{2g-2}$, where an element $(\epsilon_v)\in \Z_2^{2g-2}$ acts on $U(1)^{d(g)}$ by the formula $e^{2i\pi x_e}\mapsto (-1)^{\epsilon_v\epsilon_{v'}}e^{2i\pi x_e}$ for $v$ and $v'$, the indices of the pants bounding $C_e$.

If we choose a Lagrangian submanifold $L$ of $\ba{\boS}_g$ transverse to the fibres of the torus action and mapping diffeomorphically on $V\subset U_g^{gen}$ through $h$, then we can define canonical coordinates on $h^{-1}(V)$ by setting $x^e=0$ on $L$ and $y_e=h_e$.

The volume form is given on $h^{-1}(V)$ by  $\prod dy_e \prod dx_e$.

\end{thm}

We come back to the integral of the function associated to $\gamma$ on the moduli space $\boS_g$. Recall that $\gamma$ was adapted to the pants decomposition.  This means that $\gamma$ is the union of parallel curves $C_e$ with multiplicity $m_e$. The function $f_{\gamma}$ is then defined by $f_\gamma(\rho)=\prod_e (-\tr \rho(C_e))^{m_e}=\prod_e (-2\cos(\pi h_e(\rho)))^{m_e}=F(h)$, where $F$ is the function of Lemma \ref{calcul}.

As this function only depends on the values of $h$, we can perform the integration on its fiber first. 
The fiber is isomorphic to $U(1)^{d(g)}/\Z_2^{2g-2}$. Hence $U(1)^{d(g)}$  is a Riemannian covering over the fiber and has volume equal to 1. To find the volume of the fiber, it is then sufficient to find the degree of this covering. Let $G$ be the graph associated to the pants decomposition. The degree of the covering is equal to the dimension of the $\Z_2$-subspace of $C^1(G,\Z_2)$ 
generated by the family of vectors $u_v=e_a+e_b+e_c$ for each pant $v$ bounding circles $a,b$ and $c$. This subspace is the image of the boundary map $d:C^0(G,\Z_2)\to C^1(G,\Z_2)$. Its dimension is then complementary to the dimension of $H^1(G,\Z_2)$ which is $g$. We find that the dimension is $d(g)-g$, hence the covering has degree $2^{d(g)-g}$ and the volume of the fiber is $2^{g-d(g)}$.

We finally obtain $\int_{\boS_g} f_{\gamma} dV = 2^{g-d(g)}\int_{U_g} F(\tau)d\tau=\langle \gamma\rangle$, which completes the proof.

\end{proof}

\nocite{*}

\bibliographystyle{alpha}
\bibliography{courbe}

\textsc{Université Pierre et Marie Curie, Analyse Algébrique, Institut de Math. de Jussieu, Case 82, 4, place Jussieu, F-75252 Paris Cedex 05}\\
{\it E-mail adress:} {\tt marche@math.jussieu.fr}
\medskip

\textsc{Université Denis Diderot, Topologie et Géométrie algébriques, Institut de Math. de Jussieu, Case 7012, 2,place Jussieu, 75251 Paris Cedex 05}\\
{\it E-mail adress:} \tt{nariman@math.jussieu.fr}
\end{document}

%% file: kauf.pstex_t
\begin{picture}(0,0)%
\includegraphics{kauf.pstex}%
\end{picture}%
\setlength{\unitlength}{3947sp}%
\begingroup\makeatletter\ifx\SetFigFont\undefined%
\gdef\SetFigFont#1#2#3#4#5{%
  \reset@font\fontsize{#1}{#2pt}%
  \fontfamily{#3}\fontseries{#4}\fontshape{#5}%
  \selectfont}%
\fi\endgroup%
\begin{picture}(3358,1893)(266,-1545)
\put(1216,-141){\makebox(0,0)[lb]{\smash{\SetFigFont{8}{9.6}{\rmdefault}{\mddefault}{\updefault}{\color[rgb]{0,0,0}$=A$}%
}}}
\put(2395,-141){\makebox(0,0)[lb]{\smash{\SetFigFont{8}{9.6}{\rmdefault}{\mddefault}{\updefault}{\color[rgb]{0,0,0}$+A^{-1}$}%
}}}
\put(1135,-1150){\makebox(0,0)[lb]{\smash{\SetFigFont{8}{9.6}{\rmdefault}{\mddefault}{\updefault}{\color[rgb]{0,0,0}$=-(A^2+A^{-2})$}%
}}}
\end{picture}

%% file: bretzel.pstex_t
\begin{picture}(0,0)%
\includegraphics{bretzel.pstex}%
\end{picture}%
\setlength{\unitlength}{3947sp}%
\begingroup\makeatletter\ifx\SetFigFont\undefined%
\gdef\SetFigFont#1#2#3#4#5{%
  \reset@font\fontsize{#1}{#2pt}%
  \fontfamily{#3}\fontseries{#4}\fontshape{#5}%
  \selectfont}%
\fi\endgroup%
\begin{picture}(3134,2362)(-25,-2684)
\put(1423,-926){\makebox(0,0)[lb]{\smash{\SetFigFont{8}{9.6}{\rmdefault}{\mddefault}{\updefault}{\color[rgb]{0,0,0}$D_e$}%
}}}
\put(1390,-1743){\makebox(0,0)[lb]{\smash{\SetFigFont{8}{9.6}{\rmdefault}{\mddefault}{\updefault}{\color[rgb]{0,0,0}$D_f$}%
}}}
\put(1417,-2641){\makebox(0,0)[lb]{\smash{\SetFigFont{8}{9.6}{\rmdefault}{\mddefault}{\updefault}{\color[rgb]{0,0,0}$D_g$}%
}}}
\end{picture}